\documentclass[12pt]{article}
\usepackage{amsmath}
\usepackage{amssymb}
\usepackage{mathrsfs}
\usepackage{graphicx}

\begin{document}

\title{Explicit Examples of Strebel Differentials}
\author{Philip Tynan}

\maketitle

\begin{section}{Introduction}

In this paper, we investigate Strebel differentials, which are a special class of quadratic differentials on Riemann surfaces.  We now give a few definitions. \\

\underline{Definition}: Given a meromorphic quadratic differential $\omega$ on a Riemann surface $X$, a curve $\gamma: I \to X$ which avoids the poles of $\omega$ is part of a \textit{horizontal leaf} of $\omega$ if $f(\gamma(t)) (\gamma'(t))^2 > 0$ for each $t \in I$, where $\omega = f(z) dz^2$ in local coordinates. \\

\underline{Definition}: A meromorphic quadratic differential is a \textit{Strebel differential} if the union of its noncompact leaves forms a set of measure $0$. \\

\underline{Definition}: Let $\omega$ be a meromorphic quadratic differential on a Riemann surface $X$.  A \textit{critical trajectory} or \textit{critical horizontal leaf} of $\omega$ is a horizontal leaf of $\omega$ which tends towards a zero or pole in one direction. \\

\underline{Definition}: Let $\omega$ be a meromorphic quadratic differential.  If $a, b$ are zeroes of $\omega$ that are connected by a horizontal leaf, then a \textit{period} of $\omega$ is one of the real numbers $\displaystyle \int_a^b \sqrt \omega$, for one of the horizontal leaves connecting $a, b$. \\

Here, we will construct concrete examples of holomorphic Strebel differentials on hyperelliptic curves, as well as an algebraically defined Strebel differential on a Riemann surface defined over $\overline{\mathbb Q}$, all of whose periods are transcendental.  While the above definition of a Strebel differential is the standard definition, we need an easier criterion to check in order to show that the differentials we construct in this paper are Strebel.  The following theorem, from [5], gives us such a condition.  First, we make one more definition. \\

\underline{Definition}: The \textit{critical graph} of a meromorphic quadractic differential is the union of its zeroes, poles, and critical trajectories. \\

\underline{Theorem}: If $X$ is a Riemann surface and the genus of $X$ is greater than $1$, then a quadratic differential $\omega$ on $X$ with no poles of order $2$ or higher is Strebel iff it's critical graph is compact. \\

From [5], we also know that if $\omega$ is a Strebel differential on a Riemann surface $X$, then $\omega$ has no poles of order greater than $2$.  The property of being Strebel has some consequences that are not immediately obvious.  For example, if $\omega$ is a Strebel differential on $X$, it's critical graph gives us a unique cellular decomposition of $X$.  However, many of the theorems about Strebel differentials are non-constructive, so in this paper we will focus on constructing explicit examples. \\

We end this section with one more theorem and definition, which are cited in [3]. \\

\underline{Theorem}: A Riemann surface $X$ is algebraic iff there is a branched cover $\beta: X \to \mathbb P^1$ which is branched over only $0, 1, \infty$.  Such a map is known as a \textit{Belyi map}. \\

\end{section}

\begin{section}{Holomorphic Strebel differentials on hyperelliptic curves}
In [1], an explicit example of a Strebel differential on a genus $3$ surface was constructed.  Here, we shall construct a simpler example of a holomorphic Strebel differential on a large family of Riemann surfaces of genus $2$, which we can actually generalize to a large family of genus $g$ hyperelliptic curves, for any positive $g$.  In fact, we will see that this family is a real co-dimension $1$ subspace of the moduli space of all Riemann surfaces of genus $2$. \\

The case of genus $2$ is significant because of the following fact.  If we were to find an example of a regular Strebel differential $\omega$ on a genus $2$ surface $X$, then if $Y$ is an $n$-fold covering space of $X$, the Riemann-Hurwitz formula gives $2 - 2g(Y) = n(2 - 2g(X)) = -2n$, so $g(Y) = n+1$.  Thus, by considering $n$-fold covers of $X$, we would have an example for any positive genus. \\

From [2], we know that any smooth genus $2$ curve is isomorphic to a branched covering of $\mathbb P^1$, branched over $0, 1, \infty, \beta_1, \beta_2, \beta_3$, where $\beta_1, \beta_2, \beta_3$ are distinct elements of $\mathbb C \backslash \{0, 1\}$. \\

\underline{Theorem}: Let $X$ be a genus $2$ Riemann surface.  If $\beta_1 = \frac12 + r i$, for some real $r$, and $\beta_2, \beta_3 \in \mathbb C \backslash \{0, 1, \beta_1\}$ are distinct, then the quadratic differential form $\omega = \frac{4a^2(a^2+b^2)^2 dx^2}{((2abix-a(a+bi))^2 - a^2 (2ax - (a+bi))^2)((2abix - a(a+bi))^2 + b^2 (2ax - (a+bi))^2)}$, where $\frac14(\frac ba - \frac ab) = r$, is a holomorphic Strebel differential on $X$, expressing $X$ as as the curve $y^2 = x(x-1)(x - \frac12 - ri)(x - \beta_2)(x - \beta_3)$.  Furthermore, $\omega$ has degree two zeroes at the preimages of $\beta_2$ and $\beta_3$, and no other zeroes. \\

\underline{Proof}: Let $a, b \in \mathbb R \backslash \{0\}$, and consider the automorphism $\phi: \mathbb P^1 \to \mathbb P^1$ given by $\phi(z) = \frac{a + bi}{2a} \frac{z - a}{z - bi}$. \\

Then, \begin{eqnarray*}
g(a) & = & 0 \\
g(-a) & = & 1 \\
g(bi) & = & \infty \\
g(-bi) & = & \frac{(a+bi)^2}{4abi} = \frac12 + \frac14 (\frac ba - \frac ab) i \\
\end{eqnarray*}

Now, since the map $(\mathbb R \backslash \{0\})^2 \to \mathbb R$ given by $(a, b) \mapsto \frac ba - \frac ab$ is surjective, for any $r \in \mathbb R$ there is an automorphism of $\mathbb P^1$ taking $\{a, -a, bi, -bi\}$ to $\{0, 1, \infty, \frac12 + ri\}$, for some nonzero choice of $a, b$. \\

Thus, if $\pi: X \to \mathbb P^1$ is the branched covering which is branched at $0, 1, \infty, \frac12 + r i, \beta_2, \beta_3$, then, by choosing appropriate $a, b$, $\phi^{-1} \circ \pi: X \to \mathbb P^1$ is a branched covering, branched at $a, -a, bi, -bi, \beta_2', \beta_3'$, for some distinct $\beta_2', \beta_3' \in \mathbb P^1 \backslash \{a, -a, bi, -bi\}$.  It is not hard to check (also found in [5]) that if $\psi: Y \to Z$ is any branched covering, and $\alpha$ a meromorphic quadratic differential on $Z$, then $\alpha$ is Strebel iff $\psi^* \alpha$ is Strebel (this follows from $\psi$ being a proper, continuous map).  Therefore, it is sufficient to find a Strebel differential $\omega$ on $\mathbb P^1$, such that $(\phi^{-1} \circ \pi)^* \omega$ is holomorphic on $X$.  To find an $\omega$ such that $(\phi^{-1} \pi)^* \omega$ is holomorphic, we need only ask that $\omega$ has only simple poles, and has no poles outside of $a, -a, bi, -bi$, since $\frac{dz^2}{z}$ pulls back to $4 dz^2$ under the map $z \mapsto z^2$, and therefore second order branching over a simple pole gives a regular point for the pulled back differential.  Thus, if $\omega = \frac{dz^2}{(z^2 - a^2)(z^2 + b^2)} = f(z) dz^2$, its pullback is holomorphic. \\

It remains only to check that $\omega$ is Strebel on $\mathbb P^1$.  If $z = k \in \mathbb R$, then $f(k) = \frac{1}{(k^2 - a^2)(k^2 + b^2)}$, which is positive iff $k > a$.  Thus, if $\gamma$ is simply a path along $\mathbb R \cup \{ \infty \}$, then $\gamma'(t) \in \mathbb R$, so $\gamma'(t)^2 > 0$, therefore $f \circ \gamma(t) \gamma'(t)^2 > 0$ iff $|\gamma(t)| > a$.  And, if $z = ik, k \in \mathbb R$, then $f(ik) = \frac{1}{(-k^2 - a^2)(-k^2 + b^2)}$, which is positive iff $k > b$.  If $\gamma$ is a path along $i \mathbb R$, then $\gamma'(t) \in i \mathbb R$, so $\gamma'(t)^2 < 0$, so $f \circ \gamma(t) \gamma'(t)^2 > 0$ iff $|\gamma(t)| < b$.  Thus, we see that the critical horizontal leaves are the line going from $a$ to $-a$ which passes through infinity, and the line segment $[-bi, bi]$.  Thus, $\omega$ is Strebel. \\

Now, $\phi$ is given by the matrix $\displaystyle \left( \begin{array}{cc}
a+bi & -a(a+bi) \\
2a & -2abi
\end{array} \right)$.  Thus, $\phi^{-1}$ is given by the matrix $\displaystyle \frac{2a}{a^2 + b^2} \left( \begin{array}{cc}
-2abi & a(a+bi) \\
-2a & a+bi
\end{array} \right) \cong \left( \begin{array}{cc}
2abi & -a(a+bi) \\
2a & -(a+bi)
\end{array} \right)$. \\

Thus, \begin{eqnarray*}
(\phi^{-1})'(z) = \frac{2a(a^2+b^2)}{(2az - (a+bi))^2}
\end{eqnarray*}

So, we have: \begin{eqnarray*}
(\phi^{-1})^* \omega = & \frac{(\phi^{-1}(z))' dz^2}{((\frac{2abiz - a(a+bi)}{2az - (a+bi)})^2 - a^2)((\frac{2abiz-a(a+bi)}{2az-(a+bi)})^2 + b^2)} \\
= & \frac{4a^2(a^2+b^2)^2 dz^2}{((2abiz-a(a+bi))^2 - a^2 (2az - (a+bi))^2)((2abiz - a(a+bi))^2 + b^2 (2az - (a+bi))^2)}
\end{eqnarray*}

Pulling back under $\pi$ gives us our result. \\

We can generalize what was done in the previous section to a class of hyperelliptic curves of real codimension $1$ in the space of hyperelliptic curves.  Since a hyperelliptic curve is defined by $2g+2$ branch points, and we are free to apply an automorphism to move $3$ of them to $0, 1, \infty$ without changing the isomorphism class of our curve, then as long as one of the remaining points is moved to $\frac12 + ir$, for some real $r$, then we can pick any $2g-2$ points in the remainder of $\overline{\mathbb C}$ to get a hyperelliptic curve $\pi: C \to \mathbb P^1$.  Taking $\alpha$ to be the same as the previous $(\phi^{-1})^* \omega$ on $\mathbb P^1$, then $\pi^* \alpha$ is a holomorphic Strebel differential on $C$. \\

\end{section}

\begin{section}{Geometric interpretation of a regular Strebel differential}

From [3], we have a clear geometric interpretation for the critical graph of a Strebel differential $\omega$ with quadratic poles on a Riemann surface $X$.  The zeroes of $\omega$ are the $0$-cells, the critical trajectories the $1$-cells, and the discs that they bound are the $2$-cells, each containing a pole.  Thus, the Strebel differential gives us a cellular decomposition of $X$. \\

However, in the case of a holomorphic Strebel differential $\omega$ on a Riemann surface $X$, the picture becomes more complicated.  We can still think of the zeroes of $\omega$ as $0$-cells, and the critical leaves as $1$-cells.  But, these $1$-cells no longer bound discs, as we will now see. \\

Suppose that in local coordinates, $\omega = g(z) dz^2$.  If $\gamma: I \to X$ is a horizontal leaf, then by definition, $g(\gamma(t)) \gamma'(t)^2 > 0$.  Thus, $\sqrt{g(\gamma(t))} \gamma'(t)$ is real, and we can choose the branch so that it is always positive. \\

Now, if $\gamma(I)$ is a closed loop $C$, then $\displaystyle \int_{C} \sqrt \omega = \int_I \sqrt{g(\gamma(t))} \gamma'(t) dt > 0$.  Thus, if $C$ were homotopically trivial, we would be able to apply the residue theorem, and determine that $\omega$ has a pole.  But, since $\omega$ is holomorphic, this must mean that $C$ is homotopically nontrivial. \\

From [5], we know that in general, if $\Gamma$ is the critical graph of $\omega$, then $X \backslash \Gamma$ is a disjoint union of (possibly degenerate annuli).  The above argument shows us that in the case of a regular Strebel differential, all of the annuli are all nondegenerate. \\

We now examine the critical graphs of the Strebel differentials from the previous section.  If $\pi: X \to \mathbb P^1$ is a branched covering, and $\omega$ is a meromorphic nonvanishing Strebel differential on $\mathbb P^1$, the horizontal leaves of $\pi^* \omega$ are simply the preimages under $\pi$ of the horizontal leaves of $\pi$.  More specifically, the critical horizontal leaves of $\pi^* \omega$ are the preimages of the critical horizontal leaves of $\omega$, along with the preimages of any other horizontal leaves containing branch points.  The critical graph $\Gamma$ is then the disjoint union of these. \\

The case of critical trajectories arising from the noncritical trajectories will be split into $3$ cases in order to be properly dealt with.  The first of these is the case in which there is only a single branch point on the horizontal leaf.  This simply gives us $2$ loops that intersect \textit{non-transversely} at the single ramification point.  This is the generic situation, as most of the time, the branch points will lie on non-critical leaves, and in general will lie on different non-critical leaves. \\

Suppose now that there are $2$ branch points on one of the non-critical leaves.  The preimage of this trajectory under $\pi$ consists of $2$ loops, with $2$ \textit{transverse} intersections, one at each of the ramification points. \\

Now, suppose that there are $n \ge 3$ branch points $b_1, \ldots, b_n$ on a non-critical horizontal leaf.  The part of the graph corresponding to this consists of $n$ loops $B_1, \ldots B_n$, where $B_i$ intersects $B_{i+1}$ \textit{transversely} at the preimage of $b_i$, and the indices are taken modulo $n$. \\

Lastly, suppose that we have $n$ branch points $a_1, \ldots, a_n$ on one of the $2$ original critical leaves.  It's preimage under $\pi$ consists of a chain of $n+1$ loops $A_0, \ldots, A_n$, such that $A_{i-1}$ and $A_i$ intersect \textit{transversely} at the preimage of $a_i$, for $i = 1, \ldots, n$. \\

\end{section}

\begin{section}{Nonvanishing Strebel differentials on $\mathbb P^1$}

Let $\omega$ be a nonvanishing quadratic differential on $\mathbb P^1$.  Since $\mathcal K^{\otimes 2}$ is a degree $-4$ sheaf on $\mathbb P^1$, $\omega$ has $4$ poles (counting multiplicity).  Suppose that, in particular, $\omega$ has $4$ simple poles.  Then, $\displaystyle \omega = b \frac{dz^2}{(z-a_1)(z-a_2)(z-a_3)(z-a_4)}$, for distinct $a_i$, where $z - a_i$ is taken to be $1$ if $a_i = \infty$. \\

Now, let $\pi: E \to \mathbb P^1$ be an elliptic curve branched over $a_1, a_2, a_3, a_4$.  Then, $\pi^* \omega$ is a nonzero holomorphic differential on $E$, and is therefore equal to $\frac{dz^2}{c}$ for some $c \in \mathbb C \backslash \{0\}$. \\

Therefore, $\pi^* (\frac cb \omega)$ is Strebel, since the horizontal leaves are just the horizontal lines.  In fact, if $E$ is given by the lattice $\Gamma$, with generators $1, \tau$, then if $c'$ is some complex number, $\displaystyle \pi^* (c' \omega) = \frac{c' dz^2}{c}$ is a Strebel differential iff the slope of $\frac{c'}{c}$ is a rational multiple of the slope of $\tau$, since the horizontal leaves are simply lines with the slope of $\frac{c'}{c}$.  Thus, they will close iff the slope of $\frac{c'}{c}$ is a rational multiple of the slope of $\tau$.  In particular, this means that for most choices of $c'$, $c' \omega$ will not be Strebel. \\

Thus, for every nonvanishing differential $\omega$ with only simple poles on $\mathbb P^1$, there is a complex constant $c'$, in fact a family of such constants of real codimension $1$ in $\mathbb C$, such that $c' \omega$ is Strebel.  However, real codimension $1$ is the best we can do, as the condition above is both necessary and sufficient.  This suggests (but does not prove) that our previous result may be the best possible.  That is, it explains why it is difficult to write down an explicit example of a holomorphic Strebel differential on an arbitrary hyperelliptic curve. \\

\end{section}

\begin{section}{Constructing an algebraic Strebel differential with all transcendental periods}

In [4], a family of Strebel differentials on $\mathbb P^1$ were defined whose poles and zeroes were all at algebraic points, and whose residues were all rational (in fact, integers), but had a transcendental period.  This then gave a Strebel differential with these properties for any Riemann surface over $\overline{\mathbb Q}$.  However, the example still had algebraic periods as well.  In this section, we modify the example to construct a Strebel differential on $\mathbb P^1$ with the same conditions on its zeroes and poles, all of whose periods are transcendental, which we can pull back under a Belyi map to give us such an example on any algebraic Riemann surface. \\

Consider the differential $q_1 = - \frac{1}{\pi^2} \frac{z^2 - z + 1}{z^2(1 - z)^2}dz^2$ from [4].  It has poles of order 2 are $0, 1, \infty$ and simple zeroes are $\frac12 \pm i \frac{\sqrt3}{2}$.  Let $\Gamma$ be the graph consisting of the set of critical horizontal leaves of $q_1$.  Define 
\begin{eqnarray*}
\Gamma_1 & = & \{z \ | \ |z| = 1 \Re{z} \le \frac12 \} \\
\Gamma_2 & = & \{\frac12 + is \ | \ |s| \le \frac{\sqrt3}{2} \} \\
\Gamma_3 & = & \{ z \ | \ |z - 1| = 1, \Re{z} \ge \frac12\}
\end{eqnarray*}
Then, $\Gamma = \Gamma_1 \cup \Gamma_2 \cup \Gamma_3$, a picture of which is shown below. \\

\begin{center} \includegraphics{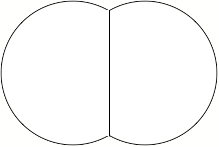} \end{center}

We saw in [4] that if $s = r \sqrt{3}$, for some rational $0 < r < \frac12$, then $\displaystyle \int_{\frac12 + is}^{\frac12 + i \frac{\sqrt3}{2}} \sqrt q$ is transcendental, so since integration of forms is invariant under pullbacks, pulling back $q_1$ over a branched cover branched only at that point and a second order pole (it is a quick check that all degree $2$ poles are preserved by any degree of branching), we obtain an algebraic Strebel differential with a transcendental period. \\

Here, we modify the construction, branching over certain algebraic points of each of the three leaves in order to produce an algebraically defined Strebel differential, all of whose periods are transcendental. \\

Consider the automorphism $h: \mathbb P^1 \to \mathbb P^1$ given by $h(z) = \frac{z - 1}{z} = 1 - \frac1z$.  We see that
$h(\frac12 + is) = \frac{-\frac12 + is}{\frac12 + is}$, so
\begin{eqnarray*}
h(\frac12 \pm i \frac{\sqrt3}{2}) & = & \frac{-\frac12 \pm i \frac{\sqrt3}{2}}{\frac12 \pm i \frac{\sqrt3}{2}} = \frac12 \pm i \frac{\sqrt3}{2} \\
h(\frac12) & = & 1 - 2 = -1 \\
h(2) & = & 1 - \frac12 = \frac12
\end{eqnarray*}

We also know that any automorphism of $\mathbb P^1$ takes circles to circles (where degenerate circles are simply lines).  Thus, since $h$ fixes two points common to all $3$ components, $h$ takes $\Gamma_2$ to $\Gamma_1$, and $\Gamma_3$ to $\Gamma_2$.  Furthermore, since $h(-1) = 2$, $h$ takes $\Gamma_1$ to $\Gamma_3$, so this automorphism has order $3$. \\

Now, suppose $c$ is a point on $\Gamma_2$.  Then, for a quadratic differential $q$ on $\mathbb P^1$, $\displaystyle \ell(c) = \int_{h(c)}^{\frac12 + i \frac{\sqrt3}{2}} \sqrt q = \int_{c}^{\frac12 + i \frac{\sqrt3}{2}} h^* \sqrt{q} = \int_{c}^{\frac12 + i \frac{\sqrt3}{2}} \sqrt{h^* q}$, (where the integrals are taken to be along $\Gamma_1, \Gamma_2$, respectively) by the definition of the pullback of a differential form. \\

But, conveniently, $\displaystyle h^* q_1 = -\frac{1}{\pi^2} \frac{(1 - \frac1z)^2 - (1 - \frac1z) + 1}{(1 - \frac1z)^2 (\frac1z)^2}(h'(z) dz)^2 \\
= -\frac{1}{\pi^2} \left( \frac{1 - \frac1z + \frac{1}{z^2}}{z^2(1 - \frac1z)^2} \right) dz^2 = -\frac{1}{\pi^2} \frac{z^2 - z + 1}{z^2(z-1)^2} dz^2 = q_1$. \\

Thus, if branching at $c$ gives a transcendental period, then branching $h(c)$ gives a transcendental period as well, as will branching at $h^{-1}(c)$. \\

Thus, if $\pi: \mathbb P^1 \to \mathbb P^1$ is a degree $4$ branched cover, branched over $c, h(c), h^{-1}(c), \infty$, with maximal ramification at $\infty$, where $c = \frac12 + i s$, and $s = r \sqrt3$, for some rational $0 < r < \frac12$, then all of the periods will be transcendental.  We now find such a map as a degree $4$ polynomial map $p: \mathbb P^1 \to \mathbb P^1$. \\

Consider the polynomial $p(z) = z^4 + c_3 z^3 + c_2 z^2 + c_1 z + c_0$.  Suppose that $p$ is ramified at $b_1, b_2, b_3$.  Then, we have that $p'(z) = 4 z^3 + 3 c_3 z^2 + 2 c_2 z + c_1 = 4 (z - b_1)(z-b_2)(z-b_3)$, and we want to impose the conditions $b_i \mapsto a_i$, where $a_2 = c, a_1 = h(c), a_3 = h^{-1}(c)$.  From here, we see that: \begin{eqnarray*}
c_3 & = & -\frac43 (b_1 + b_2 + b_3) \\
c_2 & = & 2 (b_1 b_2 + b_2 b_3 + b_3 b_1) \\
c_1 & = & -4 b_1 b_2 b_3
\end{eqnarray*}

We do not have any constraints on $c_0$.  However, we will have it depend on the ramification points of $p$, by declaring that $c_0 = b_1 b_2 b_3$.  The purpose for this will become clear later.  Let $s_i(b)$ be the $i$th symmetric polynomial in the $b_i$'s.  This gives us $3$ equations: \begin{eqnarray*}
b_i^4 - \frac43 s_1(b) b_i^3 + 2 s_2(b) b_i^2 - 4 s_3(b) b_i + s_3(b) = a_i
\end{eqnarray*}

We now homogenize these equations, as we know that the set of homogeneous equations: \begin{eqnarray*}
b_i^4 - \frac43 s_1(b) b_i^3 + 2 s_2(b) b_i^2 - 4 s_3(b) b_i + s_3(b) b_4 = a_i b_4^4
\end{eqnarray*}
have at least 64 solutions, counting multiplicity, in $\mathbb P^3$. \\

We will show that not all of these occur at the points outside of $\mathbb C^3 = D_+(b_4)$, thus showing that we have a solution in $\mathbb C^3$. \\

Setting $b_4 = 0$ gives us: \begin{eqnarray*}
b_i^4 - \frac43 s_1(b) b_i^3 + 2 s_2(b) b_i^2 - 4 s_3(b) b_i = 0
\end{eqnarray*}

It's not hard to see that $b_3 = 0 \implies b_1 = b_2 = 0$, which is impossible.  Thus, since equation $i$ is divisible by $b_i^2$, we can cancel these, leaving us with: \begin{eqnarray*}
b_1^2 - \frac43 s_1(b) b_1 + 2 s_2(b) - 4 b_2 b_3 & = & 0 \\
b_2^2 - \frac43 s_1(b) b_2 + 2 s_2(b) - 4 b_3 b_1 & = & 0 \\
b_3^2 - \frac43 s_1(b) b_3 + 2 s_2(b) - 4 b_1 b_2 & = & 0
\end{eqnarray*}

A simple check tells us that the only solution to this is $(1, 1, 1)$ (and all constant multiples of this).  Going back to our original equations, we consider the affine component $D_+(b_3) \subset \mathbb P^3$, which contains this solution. \\

Using the natural coordinates in the affine component $D_+(b_3)$, the partial derivatives at this point are $(\frac23, -\frac13, 1), (-\frac13, \frac23, 1), (-\frac13, -\frac13, 1)$, which are linearly independent over $\mathbb C^3$.  Note that this is where the importance of our choice of $c_0$ comes in.  This tells us that these three hypersurfaces intersect transversely at this point.  But, it is impossible for $3$ degree $4$ surfaces to intersect transversely at a single point.  Thus, the point $(1 : 1 : 1 : 0)$ cannot be the only solution to this system of equations, and therefore there is a solution in $\mathbb C^3  = D_+(b_4)$.  And, since these equations all have coefficients in $\mathbb Q$, the solution is contained in $\overline{\mathbb Q}^3$, so $p$ is an algebraic map $\mathbb P^1 \to \mathbb P^1$. \\

Thus, all of the preimages of $0, 1, \infty, \frac12 \pm i \frac{\sqrt3}{2}, c, h(c), h^{-1}(c)$ are algebraic, so $p^* q_1$ has all algebraic poles and zeroes, and all of its residues are integer multiples of those of $q_1$, and are therefore integers.  Furthermore, its periods are all transcendental, as shown by the computation in [4], since all are equal to $\pm L$ or $\pm (1 - L)$, where $L = \ell(\frac12 + is)$.  Thus, we have constructed an entirely algebraically defined Strebel differential $\omega$ on $\mathbb P^1$ with all transcendental periods. \\

Now, we can also apply an automorphism $\psi$ to $\mathbb P^1$ to move two of the poles that are not at $\infty$ to the points $0, 1$, and let $\omega' = \psi^* \omega$.  If $C$ is any algebraic Riemann surface, there exists a Belyi map $\beta: C \to \mathbb P^1$.  Then, $\beta^* \omega'$ is a Strebel differential on $C$, and since $\beta$ only branches over poles of $\omega$, we can now apply the rest of the argument in [4], giving us our desired differential on $C$. \\

\end{section}

\begin{section}{References}

\noindent [1] N. Ya. Amburg, \textit{An example of a regular Strebel differential}, Communications of the Moscow Mathematical Society, (2002), pp. 987-988. \\

\noindent [2] Robin Hartshorne, \textit{Algebraic Geometry}, Springer, 2006. \\

\noindent [3] Motohico Mulase and Michael Penkava, \textit{Ribbon graphs, quadratic differentials on Riemann surfaces, and algebraic curves defined over $\overline{\mathbb Q}$}, Asian Journal of Mathematics, 2 (1998), pp. 875-920. \\

\noindent [4] Motohico Mulase and Michael Penkava, \textit{Periods of Strebel differentials and algebraic curves defined over the field of algebraic numbers}, Asian Journal of Mathematics, 6 (2002), pp. 743-748. \\

\noindent [5] Kurt Strebel, \textit{Quadratic Differentials}, Springer-Verlag, 1984. \\

\end{section}

\end{document}